\documentclass{article}%
\usepackage{amsmath}
\usepackage{graphicx}%
\usepackage{amsfonts}%
\usepackage{amssymb}
\newtheorem{theorem}{Theorem}

\newtheorem{corollary}[theorem]{Corollary}

\newtheorem{definition}[theorem]{Definition}

\newtheorem{lemma}[theorem]{Lemma}

\newtheorem{proposition}[theorem]{Proposition}

\newenvironment{proof}[1][Proof]{\textbf{#1.} }{\ \rule{0.5em}{0.5em}}

\begin{document}
\bigskip

\begin{center}
{\large Group Actions on Vector Bundles }

{\large on the Projective Plane}

\bigskip

Satyajit Karnik
\end{center}

\bigskip

\textbf{Abstract}{\small : We study the action of the group of automorphisms
of the projective plane on the Maruyama scheme of sheaves }$\mathcal{M}%
_{P^{2}}(r,c_{1,}c_{2})${\small \ of rank }$r${\small \ and Chern classes
}$c_{1}=0${\small \ and }$c_{2}=n${\small \ and obtain sufficient conditions
for unstability in the sense of Mumford's geometric invariant theory. The
conditions are in terms of the splitting behaviour of sheaves when restricted
to lines in the projective plane. A very strong parallel is observed with
Mumford's theory for the action of the automorphism group of the projective
plane on the spaces of curves of a fixed degree.}

\bigskip

\textbf{1. Introduction}

The problem we study belongs mainly to the domain of geometric invariant
theory. The action of the group of automorphisms of projective space on its
hypersurfaces, or more generally on the Hilbert scheme of subvarieties with a
fixed Hilbert polynomial has been studied since Mumford. What if we replace
the Hilbert scheme of subvarieties with a suitable moduli space of sheaves or
vector bundles? The group of automorphisms of the projective space $PGL(n)$
acts naturally on sheaves by pullback. What is the nature of the semistable or
unstable points in the sense of Mumford's GIT? This question was first posed
by C.Simpson \cite{Simpson}. We study this problem for sheaves on the
projective plane.

For rank $2$ vector bundles the problem is intimately related to the action of
$SL(3)$ on the space of curves of degree $n$ in $P^{2}$. Mumford has given a
complete classification of stable, semistable and unstable points for this
action when the degree of the curve is less than $7$. In general, for any
degree, the description of the semistable and unstable points becomes
exceedingly difficult. However two fairly general statements about the `best
points' and the `worst points' remain true in all cases:

Let $C$ be a curve of degree $n$ in $P^{2}$.

1) If \ $n\geq3$, and $C$ is nonsingular, then $C$ is stable under the $SL(3)$ action,

2) If $C$ has a point of multiplicity $>$ $\frac{2n}{3}$ then $C$ is unstable.

\bigskip

The purpose of this paper is to show that there is a remarkably similar
picture when we consider the action of $SL(3)$ on sheaves on the projective
plane. As a suitable candidate for the space of sheaves on the projective
plane, we use Maruyama's scheme $\mathcal{M}_{P^{2}}(r,c_{1,}c_{2})$, of
$S$-equivalence classes of Geiseker semistable sheaves on $P^{2}$ with fixed
rank $r$ and Chern classes $c_{1}$and $c_{2}$. The precise definitions will
follow in the subsequent section. On the open subset of $\mathcal{M}_{P^{2}%
}(r,c_{1,}c_{2})$ consisting of Geiseker \textit{stable }sheaves
$S$-equivalence reduces to isomorphism of sheaves. We now state our main theorem

\begin{theorem}
Let $SL(3)$ act on $\mathcal{M}_{P^{2}}(r,0_{,}n)$ by pullback of
sheaves.Under a suitable linearization of the action the following statements
hold;\newline 1) For $n\geq3$, a generic point in $\mathcal{M}_{P^{2}}%
(r,0_{,}n)$ is stable,\newline 2) For $F\in\mathcal{M}_{P^{2}}(r,0_{,}n)$, if
there is a line $l\subseteq P^{2}$such that $F|_{l}={\large \oplus}%
\mathcal{O}(d_{i})$ with $\Sigma_{d_{i}\geq0}d_{i}>\frac{2n}{3},$then $F$ is unstable.
\end{theorem}

In particular if $F$ is rank $2$ with $c_{1}=0$ and $c_{2}=n$ such that there
is a line $l\subseteq P^{2}$ with $F|_{l}=\mathcal{O}(d)\oplus\mathcal{O}(-d)$
and $d>\frac{2n}{3}$ then $F$ is unstable under the $SL(3)$ action on
$\mathcal{M}_{P^{2}}(2,0_{,}n)$. Note the similarity with the classical result
for curves. The problem for vector bundles of rank $2$ is actually very
closely related with the case of curves. The connection is through the so
called Barth map \cite{Barth1}, that assigns to each vector bundle of rank $2$
in $\mathcal{M}_{P^{2}}(2,0_{,}n)$ a curve of degree $n$ in $P^{2\ast}$,
called its curve of jump lines. The map is $SL(3)$ equivariant and a
comparison theorem of Reichstein \cite{Reichstein} allows us to compare the
stability of bundles with the stability of curves. In higher rank, this
facility of the Barth map is lost and a new technique becomes necessary to
study stability. To do this we consider actions on suitable master spaces
obtained by monadic constructions of Beilinson \cite{Beilinson}, Barth
\cite{Barth2} and Horrocks \cite{Horrocks}. This is the main contribution of
the paper.

The paper is organised as follows. We start with recalling some basic
definitions and preliminaries in the next section. Then we work out the
special case for rank $2$ sheaves, show how it gets related to the case of
curves, and thus explain what to expect in general. In the final section we
give a proof for higher rank.

\textit{Acknowledgements}: I deeply thank my advisor Prof. Igor Dolgachev for
pointing out this problem to me and for his continuing guidance and support. I
would like to thank William Fulton, Mike Roth, Niranjan Ramachandran and Fedor
Bogomolov for valuable comments and encouragement.

\bigskip

\textbf{2}. \textbf{Definitions and Preliminaries}.

Here we recall some basic definitions and preliminaries that will be used
throughout. We begin with the fundamental notion of a semistable sheaf on a
projective space. It is a well known fact that $P_{F}(m)=$dim $H^{0}%
(P^{n},F(m))$ is a polynomial in $m$ for large enough $m$, called the Hilbert
polynomial of $F$. Let $r$ be tthe rank of $F$. The polynomial $p_{F}%
(m)=P_{F}(m)/r$ is called the reduced Hilbert polynomial of $F$. Define
$\mu(F{\normalsize )=}\frac{c_{1}(F)}{r}$.

\begin{definition}
A coherent torsion free sheaf $F$ is called \textit{semistable (stable) }if
for every coherent subsheaf $F^{^{\prime}}\subseteq F$ of rank $r^{^{\prime}}$
with $0<r^{^{\prime}}<r$, we have $p_{F^{^{\prime}}}\leq p_{F}$ ($<$), the
polynomials being ordered lexicographically, where coefficients are compared
starting with the highest degree term in the polynomial.
\end{definition}

On $P^{2}$ the reduced Hilbert polynomial looks like

\begin{center}
$p_{F}(m)=\frac{1}{2}m^{2}+(\frac{3}{2}+\mu(F))m+\frac{\chi(F)}{r}$.
\end{center}

Thus $p_{F^{^{\prime}}}\leq p_{F}$ means that $\mu(F^{^{\prime}})\leq
\mu{\normalsize (F)}$ and in case of equality we require $\chi(F^{^{\prime}%
})/r^{^{\prime}}\leq\chi(F)/r$.

The Riemann Roch theorem on $P^{2}$ reads
\[
\chi(F)=r+\frac{1}{2}c_{1}(c_{1}+3)-c_{2}%
\]
where $r$ is the rank of $F$ and $c_{1}$ and $c_{2}$ are the first and second
Chern classes respectively.

\begin{proposition}
Let $F$ be a non trivial semistable sheaf on $P^{2}$ of rank $r$ and
$c_{1}=0.$Then $H^{0}(F(-1))=H^{0}(F)=H^{2}(F(-1))=H^{2}(F)=0.$ Also
$h^{1}(F(-2))=h^{1}(F(-1))=c_{2},$ and $h^{1}(F)=c_{2}-r.$ In particular,
$c_{2}\geq r.$
\end{proposition}

\begin{proof}
The proof is a straightforward application of definitions and the Riemann roch
theorem.We first show that $H^{0}(F(-1))=0$. Suppose not. Then there is a
non-trivial map $\mathcal{O}\rightarrow F(-1).$ If $G$ is the image of
$\mathcal{O}$, then since $F(-1)$ and $\mathcal{O}$ are semistable we must
have $p_{\mathcal{O}}\leq p_{G}\leq p_{F(-1)}$.But this requires
$\mu(\mathcal{O)\leq}\mu(F(-1)),$that is $0\leq-1,$ which is false. Hence
$H^{0}(F(-1))=0$. Also, by Serre duality $H^{2}(F(-1))\ \simeq H^{0}(F^{\ast
}(-2))^{\ast}=$ Hom$(F,\mathcal{O}(-2))^{\ast}$. A non-trivial map
$F\rightarrow\mathcal{O}(-2)$ results in $\mu(F)\leq\mu(\mathcal{O}(-2))$,
that is $0\leq-2$, which is false. Hence $H^{2}(F(-1))=0$. Therefore
$\chi(F(-1))=-$dim $H^{1}(F(-1))\leq0$, and by Riemann Roch $\chi
(F(-1))=-c_{2}$. This implies that $c_{2}\geq0$ and $h^{1}(F(-1))=c_{2}$.By
exactly the same line of argument we get $\ h^{1}(F(-2))=c_{2}$.

Now we show that $H^{0}(F)=0.$ Suppose not. Then there is a non-trivial map
$\mathcal{O}\rightarrow F,$ leading to $p_{\mathcal{O}}\leq p_{F}.$ Since we
have $\mu(\mathcal{O})=\mu(F)$ semistability requires that $\chi
(\mathcal{O})\leq\chi(F)/r,$ that is $1\leq\frac{r-c_{2}}{r}=1-\frac{c_{2}}%
{r}$. Now since $c_{2}\geq0,$ we have $c_{2}=0.$But there are no semistable
sheaves other that the trivial one with $c_{1}=c_{2}=0,$ and we have assumed
that $F$ is non-trivial. Hence $H^{0}(F)=0$.Similarly $H^{2}(F)=0$ and hence
by Riemann Roch $\chi(F)=r-c_{2}=-h^{1}(F)$, so that $h^{1}(F)=c_{2}-r$.
\end{proof}

To give interesting and useful examples of stable bundles we will first look
at bundles $F$ with $c_{1}=0$ that acquire a section after a single twist by
$\mathcal{O}(1)$, that is, $H^{0}(F(1))\neq0$. Such bundles are called
Hulsbergen, after Wilfred Hulsbergen who initiated their study in his thesis
\cite{Hulsbergen}. If $s$ is a section of $F(1)$ then we get an exact
sequence
\[
0\longrightarrow\mathcal{O}\longrightarrow F(1)\longrightarrow\mathcal{I}%
_{Z}(2)\longrightarrow0\text{.}%
\]
where $\mathcal{Z}$ is the zero scheme for $s$. We will first prove the following

\begin{proposition}
The Hulsbergen bundle $F$ is Geiseker stable if and only if not all points of
$\mathcal{Z}$ \ lie on a line in $P^{2}$.
\end{proposition}

\begin{proof}
The proof rests on a simple observation that a rank $2$ bundle on $P^{2}$ with
$c_{1}=0$ is stable if and only if $H^{0}(F)=0$. If $F$ is stable then clearly
$H^{0}(F)=0$ by the above proposition. Conversely, suppose that $H^{0}(F)=0$.
Let $G\subseteq F$ be a rank $1$ subsheaf. Because $F$ is reflexive we may
assume that $G$ is in fact a bundle. Hence $G=\mathcal{O}(d)$ for some integer
$d$. But this means that $F(-d)$ has a section and hence $H^{0}(F(-d))\neq0$.
This implies that $d<0$ and hence $\mu(G)<\mu(F)$ proving stability.

Now consider the Hulsbergen bundle $F$ with the above sequence. Twisting by
$\mathcal{O}(-1)$ we see that $H^{0}(F)=0$ if and only if $H^{0}%
(\mathcal{I}_{Z}(1))=0$, which is equivalent to the geometric condition that
not all points of $Z$ lie on a line in $P^{2}$. The proof is complete.
\end{proof}

\bigskip

Let $F$ be a semistable sheaf on $P^{2}$ with reduced Hilbert polynomial
$p_{F}$.Then there is a filtration $F_{0}=(0)\subseteq F_{1}\subseteq
...\subseteq F_{n}$, such that the quotients $F_{i}/F_{i-1}$ are stable with
reduced Hilbert polynomial $p_{F}$. Define $gr(F)=\oplus F_{i}/F_{i-1}$. It
can be shown that $gr(F)$ is independent of the filtration chosen for $F$.

\begin{definition}
Two semistable sheaves $F$ and $G$ are said to be $S-$equivalent if
$gr(F)\simeq gr(G).$
\end{definition}

This is an equivalence relation on the set of semistable sheaves. On stable
sheaves $S-$ equivalence reduces to isomorphism of sheaves.

\begin{theorem}
(Maruyama \cite{Maruyama}) There exists a moduli space $\mathcal{M}_{P^{2}%
}(r,c_{1,}c_{2})$ of $S-$ equivalence classes of semistable sheaves of rank
$r$ and Chern classes $c_{1}$ and $c_{2}$. The moduli space is projective. If
$(\chi,r,c_{1})=1$ the moduli space is smooth and fine. The dimension of
$\mathcal{M}_{P^{2}}(r,c_{1,}c_{2})$ is $r^{2}(2\Delta-1)+1$, where
$\Delta=\frac{2rc_{2}-(r-1)c_{1}^{2}}{2r^{2}}$.
\end{theorem}

\begin{theorem}
(Drezet, Le Potier\cite{Lepot2}) The space $\mathcal{M}_{P^{2}}(r,c_{1,}%
c_{2})$ is irreducible.
\end{theorem}

\begin{theorem}
(Drezet \cite{Drezet}) The Picard group Pic$\mathcal{M}_{P^{2}}(r,c_{1,}%
c_{2})$ is isomorphic to $\mathbb{Z}$ if $r=n$, and is isomorphic
to\ $\mathbb{Z\oplus Z}$ if $r\neq n$.
\end{theorem}

Next we state some basic definitions and facts about group actions.

Let $G$ be a reductive algebrac group acting linearly on a vector space
$V=\mathbb{C}^{n+1}$. Then $G$ acts on the projective space $P^{n}=P(V)$. Let
$Y\subseteq P^{n}$ be a closed $G$ invariant subvariety. If such an embedding
of $Y$ in $P^{n}$ is defined by the ample line bundle $L$, then the action of
$G$ is said to be linearized with respect to $L$.

\begin{definition}
A point $a\in Y$ is said to be semi-stable if there is a homogeneous $G$
invariant polynomial of positive degree that does not vanish at $a$. Let
$Y^{ss}$ denote the set of semi-stable points of $Y$. A point $a\in Y$ is said
to be stable if it is semi-stable and if the orbit map $\varphi_{a}%
:G\rightarrow Y^{ss}$ defined by $\varphi_{a}(g)=g.a$ is proper.
\end{definition}

The following fact can be verified easily \cite{Mumford}.

\begin{proposition}
Let $a^{\prime}\in V\backslash\{0\}$ be any representative of $a\in Y$. Then
$a$ is semi-stable if $0\notin\overline{orb(a^{\prime})}$, where
$\overline{orb(a^{\prime})}$ denotes the closure of the orbit of $a^{\prime}%
$in $V$. Also, $a$ is stable if $orb(a^{\prime})$ is closed and the stabilizer
$stab(a^{\prime})$ is finite.
\end{proposition}

A non-trivial morphism of algebraic groups $\lambda:C^{\ast}\rightarrow G$ is
called a $1$-parameter subgroup of $G$. If $G$ acts linearly on $V$, then
$C^{\ast}$ acts acts on $V$ by composition. By reductivity and commutativity
of $C^{\ast}$ we get a decomposition of $V$ as $V=\oplus V_{m}$ where
$\lambda(t)(v)=t^{m}v$, for $v\in V_{m}$. Let $v\in V$ be written as $v=\Sigma
v_{m}$. Then we define
\[
\mu(v,\lambda)=-min\{m/v_{m}\neq0\}.
\]
Let $Y\subseteq P^{n}=P(V)$be a closed $G$ invariant subvariety and let $y\in
Y$. Let $y^{\prime}\in V\backslash\{0\}$ be any representative for $y$. We
put
\[
\mu(y,\lambda)=\mu(y^{\prime},\lambda)\text{.}%
\]
It is easy to check that the definition is independent of the choice of the
representative of the point $y$. The following theorems are central to all our computations.

\begin{theorem}
(Hilbert-Mumford) Let $G$ be a reductive group acting linearly on a vector
space $V=C^{n+1}$. A point $x\in P(V)$ is semi-stable (stable) under the
action of $G$ if and only if it is semi-stable (stable) under the action of
all $1$-parameter subgroups of $G$.
\end{theorem}

\begin{theorem}
(Hilbert-Mumford criterion) Let $G$ act on $V=C^{n+1}$ linearly. Let
$Y\subseteq P^{n}$ be a closed $G$ invariant subvariety. A point $y\in Y$ is
semi-stable (stable) with respect $G$ if and only if $\mu(y,\lambda)\geq0$ for
all $1$-parameter subgroups $\lambda$ of $G$.
\end{theorem}

We now recall the classical result for curves.

\begin{proposition}
Consider the action of $SL(3)$ on the space $P^{N}$ of curves of degree $n$ in
$P^{2}$. If $C\in P^{N}$ has a point of multiplicity greater than
$\frac{2n}{3}$ then $C$ is unstable.
\end{proposition}

\begin{proof}
Let $f(X,Y,Z)=0$ be the equation for the curve $C$ in $P^{2}$.Without loss of
generality we can assume that the point of multiplicity $>\frac{2n}{3}$ is
$(0;0;1)$ in the projective space $P^{N}$. Then the equation of $C$ looks like

$f(X,Y,Z)=f_{m}(X,Y,Z)+f_{m+1}(X,Y,Z)+...+f_{t}(X,Y,Z)$ where $f_{i}=\Sigma
c^{(i)}X^{a}Y^{b}Z^{c}$ with $a+b+c=n$ and $a+b=i$. In particular,
$f_{m}(X,Y,Z)=\Sigma c^{(m)}X^{a}Y^{b}Z^{c}$ with $\ a+b=m$ and not all
$c^{(m)}=0$, since then the multiplicity would not be $m$. The coordinates of
$C\in P^{N}$ are just the coefficients $c^{(j)}$ of the monomials in $X$, $Y$
and $Z$. Note that $a+b=m>\frac{2n}{3}$ and $a+b+c=n$ implies that
$c<\frac{n}{3}$.

Take the $1-$parameter group $\lambda$ that defines the action $\lambda
(t)X=tX$, $\lambda(t)Y=tY$ and $\lambda(t)Z=t^{-2}Z$ on the coordinates of the
projective plane. Under the action of this $1-$parameter group the equation
$f(X,Y,Z)$ gets transformed into

$\lambda f(X,Y,Z)=\Sigma c^{(m)}t^{a+b-2c}X^{a}Y^{b}Z^{c}+\Sigma
c^{(m+1)}t^{p+q-2r}X^{p}Y^{q}Z^{r}+...$

In all of the above sums $a+b-2c>0$, $p+q-2r>0$,...etc. That is, the power of
$t$ is always $>0$ in all the sums. This is because $a+b>\frac{2n}{3}$ and
$2c<\frac{2n}{3}$. Thus we see that if $t\rightarrow0$ the point $C$ in
$P^{N}$ tends to $0$ implying that $0$ is in the closure of the orbit of
$\tilde{C}$ a representative of $C$ in $\mathbb{C}^{N+1}$. Thus $C$ is unstable.
\end{proof}

Let $G\times H$ act on a projective variety $X$ linearized with respect to an
ample line bundle $L$, where $G$ and $H$ are reductive algebraic groups. Then
the group $H$ acts on $X$ by restriction. Let $X^{ss}(H)$ be the set of
semistable points for the action of $H$ on $X$. Let $Y=X//H$ be the GIT
quotient and let $\pi$ be the quotient map $X^{ss}(H)\rightarrow Y$. By GIT,
$Y$ comes equipped with an ample line bundle $E$ such that $\pi^{\ast
}(E)=L^{\otimes m}|_{X^{ss}(H)}$ for some $m>0$. Now $G$ acts on the quotient
space $Y$ and the action is linearized with respect to the ample line bundle
$E$. We will need to know the unstable points for this action. Let us call
this set $Y^{un}(G)$.Let $X^{un}(G\times H)$ be the locus of unstable points
for the action of $G\times H$ on $X$ linearized by the line bundle $L$.
Consider the intersection $V=X^{un}(G\times H)\cap X^{ss}(H)\subseteq
X^{ss}(H)$.

\begin{proposition}
$\pi(V)=Y^{un}(G)$.
\end{proposition}

\begin{proof}
The proof follows from the equality

$H^{0}(Y,E^{\otimes n})^{G}=H^{0}(X,L^{\otimes mn})^{G\times H}$ for all $n$,
and $m$ as described in the above paragraph. Now let $y\in Y$ be unstable. Let
$x$ be a preimage under $\pi$. By definition of an unstable point, if $s\in
H^{0}(Y,E^{\otimes n})^{G}$ then $s(y)=0$, for any $n$. But from the above
equality of groups, if $t\in H^{0}(X,L^{\otimes k})^{G\times H}$ then
$t(x)=0$, for any $k$. And this means that $x\in X^{ss}(H)$ is $G\times H$
unstable. The converse also follows.
\end{proof}

\bigskip

\textbf{3. The Barth Map and GIT}

In this section we will show how the theory of jump lines can be used to study
stability of sheaves. We begin with the classical notion of a jump line.

\begin{definition}
Let $F$ be a torsion free coherent sheaf on $P^{2}$ with $c_{1}=0$ and rank
$r$. A line $l\subseteq P^{2}$ is called a jump line for $F$ if $F|_{l}\neq
r\mathcal{O}$.
\end{definition}

The following proposition is a corollary of a famous theorem of Grauert and
Mullich\cite{Grauert}.

\begin{proposition}
Let $F$ be a semistable sheaf of rank $2$ on $P^{2}$ with $c_{1}(F)=0$. Then
for a general line $l\subseteq P^{2}$ we have $F|_{l}=\mathcal{O}%
\oplus\mathcal{O}$. That is, a generic line is non jumping.
\end{proposition}

The set of jump lines, in fact, has a better structure, as is proved by Barth.

\begin{theorem}
(W.Barth \cite{Barth2}) Let $F$ be semistable of rank $2$ on $P^{2}$ with
$c_{1}=0$. The set of jump lines of $F$, $J(F)$ $\subseteq P^{2^{\ast}}$ is a
curve of degree equal to $c_{2}$.
\end{theorem}

We explain the scheme structure of $J(F)$ now. Consider the incidence manifold
$I\subseteq P^{2}\times P^{2^{\ast}}$ defined by $I=\{(x,l)/x\in l\}$. Let
$p_{1}$ and $p_{2}$ be the projections $p_{1}:I\rightarrow P^{2}$ and
$p_{2}:I\rightarrow P^{2^{\ast}}$. Consider the sheaf $R^{1}p_{2\ast}%
p_{1}^{\ast}F(-1)$. Define $J(F)$ to be the zero scheme of the Fitting ideal
associated to the sheaf $R^{1}p_{2\ast}p_{1}^{\ast}F(-1)$. To see that $J(F)$
as defined above is indeed the set of jump lines, we must prove that a line
$\ l$ is a jump line if and only if $H^{1}(F(-1)|_{l})\neq0$.But this is clear
since a line $l$ is jumping for $F$ if and only if there is a map
$F(-1)|_{l}\rightarrow\mathcal{O}_{l}(-2)\rightarrow0$, which in turn gives a
surjection $H^{1}(F(-1)|_{l})\rightarrow H^{1}(\mathcal{O}_{l}(-2))\rightarrow
0$, and since $H^{1}(\mathcal{O}_{l}(-2))\neq0$, we have $H^{1}(F(-1)|_{l}%
)\neq0$ and we are done.

\begin{proposition}
The Barth map, b, which is the assignment $F\rightarrow J(F)$ is well defined
on $S-$ equivalence classes of semistable sheaves. Hence we have a regular map
on the moduli space b: $\mathcal{M}_{P^{2}}(2,0_{,}n)\rightarrow P^{N}$, where
$N=\frac{(n+1)(n+2)}{2}-1$.
\end{proposition}

\begin{proof}
Let $F$ be semistable. Then there is a Jordan-Holder filtration $(0)\subseteq
F_{1}\subseteq F$ such that both quotients $F_{1}$ and $F/F_{1}=G_{1}$ have
the same reduced Hilbert polynomial as $F$. In particular, $c_{1}=0$ for both
$F_{1}$ and $G_{1}$. We claim that $J(F)=J(F_{1})\cup J(G_{1})$. In view of
this let $l\in J(F_{1})\cup J(G_{1})$. Then either $H^{0}(F_{1}(-1)|_{l}%
)\neq0$ or $H^{0}(G_{1}(-1)|_{l})\neq0$ implying that $H^{0}(F(-1)|_{l})\neq
0$. Conversely, if $l\notin J(F_{1})\cup J(G_{1})$ then $H^{0}(F_{1}%
(-1)|_{l})=H^{0}(G_{1}(-1)|_{l})=H^{0}(F(-1)|_{l})=0$ which means that $l$ is
not a jump line for $F.$ Thus if $f=0$ and $g=0$ are the equations of
$J(F_{1})$ and $J(G_{1})$ respectively then $J(F)$ has equation $fg=0$. Since
the Jordan-Holder filtration is unique upto isomorphism of quotients, $J(F)$
is well defined on $S-$ equivalence classes of sheaves. this completes the proof.
\end{proof}

The Hulsbergen bundles offer examples of curves of jump lines. Let $F$ be a
Hulsbergen bundle defined by the sequence
\[
0\longrightarrow\mathcal{O}\longrightarrow F(1)\longrightarrow\mathcal{I}%
_{Z}(2)\longrightarrow0
\]
where $\mathcal{Z}$ is a zero scheme in $P^{2}$. Suppose that $d+1$ point of
$\mathcal{Z}$ lie on a line $l\subseteq P^{2}$.When restricted to this line
$l$ the sheaf $\mathcal{I}_{Z}$ is simply $\mathcal{O}_{l}(-d-1)$. The above
defining sequence gives a surjection $F|_{l}\rightarrow\mathcal{O}%
(-d)\rightarrow0$. Since $F$ is rank $2$ and $c_{1}(F)=0$, the kernel of the
map must be $\mathcal{O}(d)$. Since $H^{1}(P^{1},\mathcal{O}(2d))=0$ for
$d\geq0$, the restricted sequence splits and we get the decomposition
$F|_{l}=\mathcal{O}(d)\oplus\mathcal{O}(-d)$. Hence if $d\geq1$ the line $l$
is a jump line for $F$.

One can in fact write down equations for the curve of jump lines. Let $Z$ be
the zero scheme of $n+1$ distinct points $\{x_{1},...,x_{n+1}\}$ in $P^{2}$ no
three of which are collinear. Let $L_{1},...,L_{n+1}$ be linear forms that are
equations for the lines $x_{1}^{\ast},...,x_{n+1}^{\ast}$ in $\ P^{2\ast}$.
Let $f_{i}=\Pi_{j\neq i}L_{j}$. Let $W=span\{f_{1},...,f_{n+1}\}\subseteq
H^{0}(P^{2},\mathcal{O}(n))$. Hulsbergen proves that thee is an isomorphism
$Ext(I_{Z}(2),\mathcal{O})\simeq W=\mathbb{C}^{n+1}$ such that if $F\in
Ext(I_{Z}(2),\mathcal{O})$ corresponds to $(a_{1},...,a_{n+1})\in W$ the the
curve of jump lines of $F$, $J(F)$, is given by the equation
\[
a_{1}f_{1}+...+a_{n+1}f_{n+1}=0\text{.}%
\]
For a generic $F$, the curve of jump lines is clearly nonsingular.

In order to exploit the Barth map to understand stability of sheaves we need
the following comparison theorem of Reichstein\cite{Reichstein} on a relative
Hilbert-Mumford criterion.

\begin{theorem}
(Reichstein) Let $X$ and $Y$ be projective varieties on which a reductive
group acts. Let $f:X\rightarrow Y$ be a $G-$equivariant map. Let the actions
of $G$ on $X$ and $Y$ be linearized by ample line bundles $L$ and $M$
respectively. Let $L_{d}=f^{\ast}M^{\otimes d}\otimes L$. For large $d$ the
following holds\newline (i) If $f(x)=y$ and $y$ is unstable, then $x$ is
unstable under the linearization $L_{d}$\-$,$\newline (ii) If $f(x)=y$ and $y$
is stable, then $x$ is stable under the linearization $L_{d}$.
\end{theorem}

We now prove our first theorem for rank $2$ sheaves.

\begin{theorem}
Let $F\in\mathcal{M}_{P^{2}}(2,0_{,}n)$. Suppose there is a line $l\subseteq
P^{2}$ such that $F|_{l}=\mathcal{O}(d)\oplus\mathcal{O}(-d)$ with
$d>\frac{2n}{3}$, then $F$ is unstable for the $SL(3)$ action suitably
linearized. Also, a generic point in $\mathcal{M}_{P^{2}}(2,0_{,}n)$ is stable.
\end{theorem}

\begin{proof}
With the previous background the proof is quite easy. Consider the Barth map
$b:\mathcal{M}_{P^{2}}(2,0_{,}n)\rightarrow P^{N}$. On $J(F)$ the line $l$
corresponds to a point of multiplicity $m$, with $m\geq d$. this can be seen
in Barth's fundamental paper \cite{Barth2}. Now since $d>\frac{2n}{3}$, we
have $m>\frac{2n}{3}$. Hence $J(F)$ is unstable. Hence by Reichstein's theorem
$F$ is unstable under suitable linearization.

To prove that a generic point is stable we simply produce a stable point.
Genericness follows from the irreducibility of the moduli space. Consider a
zero scheme $\mathcal{Z}$ of length $n+1$ in $P^{2}$ such that no three points
are on a line. Let $F$ be a vector bundle on $P^{2}$ given by the Hulsbergen
sequence%
\[
0\longrightarrow\mathcal{O}\longrightarrow F(1)\longrightarrow\mathcal{I}%
_{\mathcal{Z}}(2)\longrightarrow0\text{.}
\]

Then the curve of jump lines $J(F)$ is nonsingular and hence stable by
Mumford's result. But then by Reichstein's theorem $F$ is also stable. The
proof is complete.
\end{proof}

\bigskip

\textbf{3. Higher Rank: Group Action on }$\mathcal{M}_{P^{2}}(r,0,n)$

When the rank $r\geq3$ the facility of the Barth map is lost.
Hulek\cite{Hulek2} has tried to recover the theory of jump lines for higher
rank. However, it turns out that in the event that a curve of jump lines can
be assigned to a vector bundle, a lot of information is lost in the process.
For one, the degree of the curve is much larger than $\ c_{2}$ and the
multiplicities of points on the curve are very weakly related to the splitting
of the bundle over the line. Thus in order to get sharp results we need to
adopt a different approach.

In 1970 Horrocks realized a remarkable philosophy of recovering a vector
bundle from its cohomology groups and certain maps between them. This is his
theory of monads. Monadic representations allow us to view the moduli space as
a GIT quotient of subvarieties of products of Grassmanians. In this sense it
replaces the so-called Quot scheme. This structure is more convenient to study
Picard groups and group actions. The GIT is reduced to calculations on
products of grassmanians. We start by stating a fundamental result of
Beilinson on spectral sequences.

\begin{theorem}
(Beilinson \cite{Beilinson}) Let $\ E$ be a rank $r$ torsion free sheaf on
$P^{n}$. Then there is a spectral sequence $E_{r}^{pq}$, for $0\leq q\leq n$
and $-n\leq p\leq0$, with $E_{1}$ term given by $E_{1}^{pq}=H^{q}%
(P^{n},E(p))\otimes\Omega^{-p}(-p),$which converges to
\begin{align*}
E^{i}  &  =E\text{ if }i=0,\\
E^{i}  &  =0\text{ otherwise.}%
\end{align*}
\end{theorem}

\begin{definition}
A monad is a complex of sheaves%
\[
0\longrightarrow A\longrightarrow B\longrightarrow C\longrightarrow0
\]
exact at $A$ and $C$ but not necessarily exact in the middle.
\end{definition}

\begin{proposition}
Let $E$ be a semistable sheaf of rank $r$ on $P^{2}$ with Chern classes
$c_{1}$ and $c_{2}$. Assume that $E$ is normalized, that is, $-r<c_{1}\leq0$.
Then $E$ is the middle cohomology of the monad%
\[
0\longrightarrow H^{1}(E(-2))\otimes\mathcal{O}(-1)\longrightarrow
H^{1}(E(-1))\otimes\Omega(1)\longrightarrow H^{1}(E)\otimes\mathcal{O}%
\longrightarrow0\text{.}
\]
\end{proposition}

\begin{proof}
Consider the $E_{1}$ term
\[
E_{1}^{pq}=H^{q}(P^{n},E(p))\otimes\Omega^{-p}(-p)\text{.}
\]

Since our base space is $P^{2}$, $E_{1}^{p,q}=0$ if $q>2$. Thus $q=0,1,2$ are
the only values of $q$ yielding non zero $E_{1}$terms. Moreover, by
semistability of $E$ we have $H^{0}(P^{2},E(p))=H^{2}(P^{2},E(p))=0$ for
$p=0,1,2$. Thus $E_{1}^{p,0}=E_{2}^{p,2}=0$ for $p=0,-1$ and $-2$. The
differential maps at level $1$ are defined as $d_{1}^{p,q}:E_{1}%
^{p,q}\rightarrow E_{1}^{p+1,q}$. We now have $\ker d_{1}^{-2,1}=E_{2}%
^{-2,1}=A$(say), $\ker d_{1}^{-1,1}/imaged_{1}^{-2,1}=E_{2}^{-1,1}=B$(say) and
$co\ker d_{1}^{-1,1}=E_{2}^{0,1}=C$(say). The differentials at level $2$ are
defined $d_{2}^{p,q}:E_{2}^{p,q}\rightarrow E_{2}^{p+1,q-1}$ and hence are
degenerate. The spectral sequence therefore degenerates at the $E_{2}$ term.
By Beilinsons' theorem $E_{2}^{-2,1}=E_{2}^{0,1}=0$, and $E_{2}^{-1,1}=E$\ $.
$This simply means that we have a complex
\[
0\longrightarrow E_{1}^{-2,1}\longrightarrow E_{1}^{-1,1}\longrightarrow
E_{1}^{0,1}\longrightarrow0\text{,}
\]

whose middle cohomology is $E$. Reading the terms of the complex we get the
required monad whose middle cohomology is $E$.
\end{proof}

\begin{corollary}
\ Let $E$ be semistable of rank $r$ on $P^{2}$ with $c_{1}=0$ and c$_{2}=r$.
Then $E$ occurs as the cokernel in the exact sequence
\[
0\longrightarrow\mathbb{C}^{r}\otimes\mathcal{O}(-1)\longrightarrow
\mathbb{C}^{r}\otimes\Omega(1)\longrightarrow E\longrightarrow0\text{.}%
\]
\end{corollary}

\begin{proof}
First note that dim$H^{1}(E)=-\chi(E)=c_{2}-r=0$. Adding this information in
the statement of the above proposition we see that $E$ is the cokernel in the
exact sequence of the form%
\[
0\longrightarrow H^{1}(E(-2))\otimes\mathcal{O}(-1)\longrightarrow
H^{1}(E(-1))\otimes\Omega(1)\longrightarrow E\longrightarrow0\text{.}
\]

By Riemann Roch theorem each of $H^{1}(E(-2))$ and $H^{1}(E(-1))$ are $r$
dimensional vector spaces. The result follows.
\end{proof}

\bigskip

The corollary of Beilinson's theorem helps us to construct the moduli space
$\mathcal{M}(r,0,n)$ as a GIT quotient of a certain `space of monads'. We
describe this construction below.

Let $K,H,L$ and $V$ \ be fixed vector spaces of dimensions $n,n,n-r$ and $3$
respectively. The projective plane under consideration will be $P^{2}=P(V)$.
Consider the product of grassmannians $G(n,H\otimes V)\times G(H\otimes
V^{\ast},n-r)$, where $G(n,H\otimes V)$ is the grassmannian of subspaces of
$H\otimes V$ of dimension $n$ and $G(H\otimes V^{\ast},n-r)$ is the
grassmannian of quotients of $H\otimes V^{\ast}$ of dimension $n-r$. A pair
$(K,L)$ in\ \ \ $G(n,H\otimes V)\times G(H\otimes V^{\ast},n-r)$ gives maps
$a:K\otimes\mathcal{O}(-1)\rightarrow H\otimes\Omega(1)$, and $b:H\otimes
\Omega(1)\rightarrow L\otimes\mathcal{O}$. Consider the subspace
$\mathfrak{M}\subseteq G(n,H\otimes V)\times G(H\otimes V^{\ast},n-r)$
comsisting of pairs $(K,L)$ such that $b\circ a=0$. This is the master space
of monads. The group $SL(H)$ acts on $\mathfrak{M}$. The following theorem of
Le Potier is very crucial for us;

\begin{theorem}
(Le Potier) For the polarization $(rm-n,n)$ on $G(n,H\otimes V)\times
G(H\otimes V^{\ast},n-r)$ restricted to $\mathfrak{M}$, with $m$ sufficiently
large, we have the isomorphism%
\[
\mathfrak{M}//SL(H)\simeq\mathcal{M}_{P^{2}}(r,0,n)\text{.}%
\]
\end{theorem}

The GIT quotient $\mathcal{M}_{P^{2}}(r,0,n)$ comes equipped with an ample
line bundle which we denote by $E_{m}$. Our discussion of stability is with
respect to the action of $SL(3)$ linearized with respect to $E_{m}$.

We will need the following lemma.

\begin{lemma}
Let $(K,L)$ be a pair defining a semistable sheaf $F$. Let $v\in V$ and let
$v^{\bot}\subseteq V^{\ast}$ be the annihilating space for $v$. The injection
$K\rightarrow H\otimes V$ gives a natural map $\alpha_{v}:K\rightarrow$
Hom$(v^{\bot},H)$ and this map is injective.
\end{lemma}

\begin{proof}
Suppose that there is a $v\in V$, $v\neq0$ such that the map $\alpha_{v}$ is
not injective. Then there is a $k\in K$, such that $k\neq0$ and $\alpha
_{v}(k)=0$.

Blow up $P^{2}$ at the point $x=\mathbb{C}.v\in P^{2}$ and denote the $\sigma$
process at $x$ by $\sigma:\widehat{P^{2}}\rightarrow P^{2}$. Embed
$\widehat{P^{2}}$ in $P^{2}\times P^{1}$ as an incidence manifold. Let $p$ and
$q$ denote the projections to $P^{2}$ and $P^{1}$ respectively. Let
$E\subset\widehat{P^{2}}$ be the exceptional curve. Then we have $\sigma
^{\ast}(\mathcal{O}_{P^{2}}(1))=q^{\ast}(\mathcal{O}_{P^{1}}(1))\otimes
\mathcal{L}(E)$. Now we have the exact sequence
\[
0\longrightarrow p^{\ast}F(-2)\otimes q^{\ast}(\mathcal{O}_{P^{1}%
}(-1))\longrightarrow p^{\ast}F(-1)\otimes q^{\ast}\mathcal{O}_{P^{1}%
}\longrightarrow p^{\ast}F(-1)\otimes q^{\ast}\mathcal{O}_{P^{1}}%
|_{\widehat{P^{2}}}\longrightarrow0\text{.}%
\]
The direct image of this sequence under $q$ gives the following sequence on
$P^{1}$,%
\[
0\longrightarrow q_{\ast}(p^{\ast}F(-2)\otimes q^{\ast}(\mathcal{O}_{P^{1}%
}(-1)))|_{_{\widehat{P^{2}}}}\longrightarrow K\otimes\mathcal{O}_{P^{1}%
}(-1)\longrightarrow H\otimes\mathcal{O}_{P^{1}}\text{.}%
\]
A point $l\in P^{1}$ corresponds to a line $L\subset P^{2}$ such that the
equation for $L$ is $z=0$ for $z\in v^{\bot}\subseteq V^{\ast}$. Thus the rank
one subsheaf $k\otimes q_{\ast}(p^{\ast}F(-1)\otimes q^{\ast}\mathcal{O}%
_{P^{1}})$ lies in the kernel of the right hand arrow. Hence $q_{\ast}%
(p^{\ast}F(-1)\otimes q^{\ast}\mathcal{O}_{P^{1}})$ contains $\mathcal{O}%
_{P^{1}}(-1)$ as a subsheaf and we have
\[
p^{\ast}F(-1)\otimes q^{\ast}\mathcal{O}_{P^{1}}(1)\simeq\sigma^{\ast}%
F\otimes\mathcal{L}(E)^{-1}\subset\sigma^{\ast}F
\]
which implies that $\sigma^{\ast}F$ contains a nontrivial section. But this is
a contradiction since we know that by semistability $H^{0}(F)=0$. This
concludes the proof.
\end{proof}

In a similar vein we have the following proposition.

\begin{proposition}
If $V^{\prime\ast}\subseteq V^{\ast}$ is a one dimensional subspace of
$V^{\ast}$ then the image of the restricted map $H\otimes V^{\prime\ast
}\rightarrow L$ is $(0)$.
\end{proposition}

\begin{proof}
We know that $0\rightarrow L^{\ast}\rightarrow H^{\ast}\otimes V$. If $v\in
V$, we will show that $L^{\ast}\cap(H^{\ast}\otimes(v))=(0)$. The proof is
along the same lines as the previous lemma except that we need to use
properties of $F^{\ast}$ namely that it is torsion free and $H^{0}(F^{\ast
}(-1))=0$. By Serre duality $L^{\ast}=H^{1}(F)^{\ast}\simeq H^{1}(F^{\ast
}(-3))$, and similarly $H^{\ast}=H^{1}(F(-1))^{\ast}\simeq H^{1}(F^{\ast
}(-2))$. Suppose that we have $\beta\in L^{\ast}\cap(H^{\ast}\otimes(v))$,
$\beta\neq0$. Let $v^{\bot}=ann(v)\subseteq V^{\ast}$ be the annihilating
subspace in $V^{\ast}$.

Blow up $P^{2}$ at the point $x=\mathbb{C}.v\in P^{2}$, and denote this
$\sigma$ process at $x$ by $\sigma;\widehat{P^{2}}\rightarrow P^{2}$. Embed
$\widehat{P^{2}}$ in $P^{2}\times P^{1}$ as an incidence manifold. Let $p$ and
$q$ denote the projections to $P^{2}$ and $P^{1}$ respectively. Let
$E\subset\widehat{P^{2}}$ be the exceptional curve. Then we have $\sigma
^{\ast}(\mathcal{O}_{P^{1}}(1))=q^{\ast}(\mathcal{O}_{P^{1}}(1))\otimes
\mathcal{L}(E)$. Now we have the exact sequence
\[
0\longrightarrow p^{\ast}F^{\ast}(-3)\otimes q^{\ast}(\mathcal{O}_{P^{1}%
}(1))\longrightarrow p^{\ast}F^{\ast}(-2)\otimes q^{\ast}\mathcal{O}_{P^{1}%
}\longrightarrow p^{\ast}F^{\ast}(-2)\otimes q^{\ast}\mathcal{O}_{P^{1}%
}|_{\widehat{P^{2}}}\longrightarrow0\text{.}%
\]
The direct image under $q$ gives the following sequence on $P^{1}$,
\[
0\longrightarrow q_{\ast}(p^{\ast}F^{\ast}(-2)\otimes q^{\ast}\mathcal{O}%
_{P^{1}})|_{\widehat{P^{2}}}\longrightarrow L^{\ast}\otimes\mathcal{O}_{P^{1}%
}(-1)\longrightarrow H^{\ast}\otimes\mathcal{O}_{P^{1}}\text{.}%
\]
A point $l\in P^{1}$ corresponds to a line $L\subseteq P^{2}$ such that the
equation for $L$ is $z\in v^{\perp}\subseteq V^{\ast}$.

Thus the rank one subsheaf $k\otimes q_{\ast}(p^{\ast}F^{\ast}(-2)\otimes
q^{\ast}\mathcal{O}_{P^{1}})$ lies in the kernel of the rightmost arrow.Thus
$q_{\ast}(p^{\ast}F^{\ast}(-2)\otimes q^{\ast}\mathcal{O}_{P^{1}})$ containes
$\mathcal{O}_{P^{1}}(-1)$ as a subsheaf and hence
\[
p^{\ast}F^{\ast}(-2)\otimes q^{\ast}F^{\ast}(1)\simeq\sigma^{\ast}F^{\ast
}(-1)\otimes L(E)^{-1}\subset\sigma^{\ast}F^{\ast}(-1)\text{.}%
\]
This implies that $\sigma^{\ast}F^{\ast}(-1)$ contains a non trivial section.
But this is a contradiction since we know that $H^{0}(F^{\ast}(-1))=0$. This
concludes the proof.
\end{proof}

\bigskip

We now begin to study actions on products on grassmannians. For convenience we
will consider grassmannians of subspaces only rather than grassmannians of
quotients. Results for the latter can be obtained simply by dualizing the actions.

\begin{proposition}
Let $G=SL(H)\times SL(V)$ act on the space $Y=G(n,H\otimes V)\times
G(m,H\otimes V)$ linearized by the line bundle $(p,q)\in Pic(X)$. Then $SL(V)$
acts on $Y$ by restriction. A point $(K,L)\in Y$ is unstable with respect to
the $SL(V)$ action if and only if there is a subspace $V^{\prime}\subseteq V$
such that if $K^{\prime}=K\cap H\otimes V^{\prime}$ and $L^{\prime}=L\cap
H\otimes V^{\prime}$ then%
\[
\frac{p\text{ dim}K^{\prime}+q\text{ dim}L^{\prime}}{\text{dim}V^{\prime}%
}>\frac{p\text{ dim}K+q\text{ dim}L}{\text{dimV}}\text{.}%
\]
\end{proposition}

\begin{proof}
First assume that the pair $(K,L)$ satisfies%
\[
\frac{p\text{ dim}K^{\prime}+q\text{ dim}L^{\prime}}{\text{dim}V^{\prime}}%
\leq\frac{p\text{ dim}K+q\text{ dim}L}{\text{dimV}}%
\]

for every $V^{\prime}\subseteq V$. We will prove that $(K,L)$ is semistable.
Let $\lambda$ be a $1-$parameter subgroup of $SL(V)$. Choose subspaces
$V_{j}\subseteq V$ such that $t(v_{j})=t^{u_{j}}v_{j}$ for $v_{j}\in V_{j}$
and $t\in\lambda$. Let $F_{i}=\oplus_{p=1}^{i}H\otimes V_{p}$. Without loss of
generality we can assume that $u_{1}>u_{2}>...u_{s}$. Then we have the flag%
\[
F_{1}\subseteq F_{2}\subseteq...F_{s}=H\otimes V\text{.}
\]

Define $K_{i}=K\cap F_{i}$ and $L_{i}=L\cap F_{i}$ for $i=1,...s$. Now we know
that%
\[
-\mu((K,L),\lambda)=-p\mu(K,\lambda)-q\mu(L,\lambda)\text{.}
\]

Hence%
\[
-\mu((K,L),\lambda)=p\Sigma_{i=1}^{s}u_{i}(\text{dim}K_{i}-\text{dim}%
K_{i-1})+q\Sigma_{i=1}^{s}u_{i}(\text{dim}L_{i}-\text{dim}L_{i-1})\text{.}%
\]

The assumed inequality for $(K,L)$ gives us%
\[
(u_{i=1}-u_{i})(p\text{ dim}K_{i-1}+q\text{ dim}L_{i-1})\leq c(\text{dim}%
F_{i-1})(u_{i-1}-u_{i})
\]

for $i=2,...s$ and%
\[
u_{s}(p\text{ dim}K_{s}+q\text{ dim}L_{s})\leq c\text{ dim}F_{s}u_{s}\text{.}
\]

Adding all these inequalities we get $\mu\geq0$. By the Hilbert-Mumford
criterion this means that $(K,L)$ is a semistable point.

Conversely, suppose that there is a subspace $V^{\prime}\subseteq V$ such that
we have the inequality%
\[
\frac{p\text{ dim}K^{\prime}+q\text{ dimL'}}{\text{dim}V^{\prime}%
}>\frac{p\text{ dim}K+q\text{ dim}L}{\text{dimV}}\text{.}
\]

Let $V=V^{\prime}\oplus V^{\prime\prime}$. Take the $1-$parameter group
$\lambda\subseteq SL(V)$ that has the action defined by $t(v^{\prime
})=t^{\text{dim}V^{\prime\prime}}v^{\prime}$ for $v^{\prime}\in V^{\prime}$
and $t(v^{\prime\prime})=t^{-\text{dim}V^{\prime}}v^{\prime\prime}$ for
$v^{\prime\prime}\in V^{\prime\prime}$. Calculating $\mu$ we find that
$\mu((K,L),\lambda)<0$, which proves unstability.
\end{proof}

\bigskip

For any subspace $V^{\prime}\subseteq V$, there is an induced map
$L\rightarrow L^{\prime}\rightarrow0$, where $L^{\prime}$ is defined as
follows: The injection $V^{\prime}\rightarrow V$ gives a surjection $V^{\ast
}\rightarrow V^{\prime\ast}$ with a kernel $W^{\ast}$. We get a map $H\otimes
W^{\ast}\rightarrow L$ by restriction. If $L^{\prime\prime}$ is the image then
$L^{\prime}$ is defined as $L/L^{\prime\prime}$. We use this consturction in
the proposition below.

\begin{corollary}
Let $SL(H)\times SL(V)$ act on the space $X=G(k,H\otimes V)\times G(H\otimes
V^{\ast},m)$ where the action of $SL(V)$ on $V^{\ast}$ is defined by
transpose. Let the action be linearized by the ample line bundle $(p,q)$ on
$X$. A point $(K,L)$ is unstable for the $SL(V)$ action if and only if there
is a subspace $V^{\prime}\subseteq V$ such that if $K^{\prime}=K\cap H\otimes
V^{\prime}$and $L^{\prime}$is as defined above we have
\[
\frac{p\text{ dim}K^{\prime}+q\text{ dim }L^{\prime}}{\text{dim}V^{\prime}%
}>\frac{p\text{ dim}K+q\text{ dim}L}{\text{dim}V}\text{.}%
\]
\end{corollary}

\begin{proof}
By dualizing we can think of $G(H\otimes V^{\ast},m)$ as the grassmannian of
subspaces $G(m,H^{\ast}\otimes V)$ . Under this identification let the pair
$(K,L)$ correspond to the pair $(K,W)$. Then from the proof of the above
proposition, $(K,W)\in G(m,H^{\ast}\otimes V)$ is $SL(V)$ unstable if and only
if there is a subspace $V^{\prime}\subseteq V$ such that if $K^{\prime}=K\cap
H\otimes V^{\prime}$ and $W^{\prime}=W\cap H^{\ast}\otimes V^{\prime}$, we
have%
\[
\frac{p\text{ dim}K^{\prime}+q\text{ dim}W^{\prime}}{\text{dim}V^{\prime}%
}>\frac{p\text{ dim}K+q\text{ dim}W}{\text{dim}V}\text{.}%
\]
The vector spaces $L$ and $W$ are related by the th equality $W=L^{\ast}$. The
subspace $W^{\prime}$ of $W$ thus corresponds to a subspace $L^{\prime\ast}$
of $L^{\ast}$. Hence we get the quotient $L\rightarrow L^{\prime}$. Since
$L^{\prime}$ and $W^{\prime}$ are of the same dimension the above inequality
reduces to
\[
\frac{p\text{ dim}K^{\prime}+q\text{ dim}L^{\prime}}{\text{dim}V^{\prime}%
}>\frac{p\text{ dim}K+q\text{ dim}L}{\text{dim}V^{\prime}}%
\]
which is the stated criterion. The proof is complete.
\end{proof}

\begin{corollary}
Let $F\in\mathcal{M}_{P^{2}}(r,0,n)$ and let $(K,L)\in G(n,H\otimes V)\times
G(H\otimes V^{\ast},n-r)$ be a pair corresponding to $F$. If $(K,L)$ is such
that there is a vector subspace $V^{\prime}\subseteq V$ with
\[
\text{dim}K^{\prime}>\frac{\text{dim}V^{\prime}.\text{dim}K}{\text{dim}%
V}\text{.}%
\]
where $K^{\prime}=H\otimes V^{\prime}\cap K$. Then $(K,L)$ is $SL(V)$
unstable. Moreover, for the above inequality to hold we need dim$V^{\prime}=2$.
\end{corollary}

\begin{proof}
Note that for the inequality in the hypothesis to hold we must have
dim$V^{\prime}=2.$For if dim$V^{\prime}=1$, then dim$K^{\prime}=0$, by the
lemma proved before and the inequality does not hold. We know that $(K,L)$ is
unstable if and only if there is a subspace $V^{\prime}\subseteq V$ such that
if $K^{\prime}=K\cap H\otimes V^{\prime}$ and $L^{\prime}$ as defined before
we have
\[
\frac{(rm-n)\text{ dim}K^{\prime}+n\text{ dim}L^{\prime}}{\text{dim}V^{\prime
}}>\frac{(rm-n)\text{ dim}K+n\text{ dim}L}{\text{dim}V}%
\]
where $m>>0$. Simplifying the above inequality gives
\[
\frac{\text{dim}K^{\prime}}{\text{dim}V^{\prime}}>\frac{\text{dim}%
K}{\text{dim}V}+\frac{1}{rm-n}\left(  \frac{-n\text{ dim}L^{\prime}%
}{\text{dim}V^{\prime}}+\frac{n\text{ dim}L}{\text{dim}V}\right)  \text{.}%
\]
Since $m$ is very large the above inequality is equivalent to
\[
\text{dim}K^{\prime}>\frac{\text{dim}V^{\prime}.\text{dim}K}{\text{dim}V}%
\pm\in
\]
where $\in$ $>0$, and is very small. Since dim$K^{\prime}$ is an integer, the
inequality in the hypothesis,
\[
\text{dim}K^{\prime}>\frac{\text{dim}V^{\prime}.\text{dim}K}{\text{dim}%
V}\text{,}%
\]
implies the inequality for unstability. The proof is complete.
\end{proof}

\bigskip

This will be crucial in the final steps of the proof. Roughly, the
'undemocratic' polarization leads the second factor to eventually drop out of
the analysis.

The following theorem relates the geometry of $F$ to its $SL(V)$ unstability.

\begin{theorem}
Let $F\in\mathcal{M}_{P^{2}}(r,0,n).$ If there is a line $l\subseteq P^{2}$
such that $F|_{l}=\oplus\mathcal{O}(d_{i})\oplus_{j=1,...,p}C_{j}$ with
\[
\Sigma_{d_{i}\geq0}d_{i}+p>\frac{2n}{3}\text{.}%
\]
then $F$ is unstable under the $SL(3)$ action on $\mathcal{M}_{P^{2}}(r,0,n)$
with respect to the polarization $E_{m}$.
\end{theorem}

\begin{proof}
Assume that there is a line $l\subseteq P^{2}$with the given property of $F$.
Let $(K,L)$ be a point in the space of monads that defines $F$. Let
$l=P(V^{\prime})$, where dim $V^{\prime}=2$. We prove that $V^{\prime}$ is a
destabilizing subspace for the pair $(K,L)$. We have the exact sequence%
\[
0\longrightarrow F(-2)\longrightarrow F(-1)\longrightarrow F(-1)|_{l}%
\longrightarrow0\text{,}%
\]
from which we get the cohomology sequence%
\[
0\longrightarrow H^{0}(F(-1)|_{l})\longrightarrow K\longrightarrow
H\longrightarrow H^{1}(F(-1)|_{l})\longrightarrow0\text{.}%
\]
Then $K^{\prime}=K\cap H\otimes V^{\prime}$ is also the image of
$H^{0}(F(-1)|_{l})$ in $K$. To check this first recall that $K$ is a subspace
of $H\otimes V$. For $l\in V^{\ast}$ the map $\phi_{l}:K\rightarrow H$ is
obtained as follows. Let $k\in K$ be written as $k=\Sigma h_{i}\otimes v_{i}%
$.Then $\phi_{l}(k)=\Sigma\phi_{l}(v_{i})h_{i}$. Since $l=P(V^{\prime})$,
$\phi_{l}$ annihilates $V^{\prime}$ and hence if $k\in K\cap H\otimes
V^{\prime}$ then $\phi_{l}(k)=0$ which means that $k\in H^{0}(F(-1)|_{l})$.

Conversely, assume that $\phi_{l}(k)=0$. Think of $k$ as $k=h_{1}\otimes
v_{1}+h_{2}\otimes v_{2}+h_{3}\otimes v_{3}$ where $v_{1},v_{2}\in V^{\prime}$
and $v_{3}\notin V^{\prime}$. Then $\phi_{l}(k)=\phi_{l}(v_{3})h_{3}=0$, which
implies that $h_{3}=0$, and hence $k\in K\cap H\otimes V^{\prime}$. Thus
$K^{\prime}$ is the image of $H^{0}(F(-1)|_{l})$ in $K$. A simple computation
on $P^{1}$yeilds dim $K^{\prime}=\Sigma_{d_{i}\geq0}d_{i}+p$. Since dim $V=3$
and dim$K=n$, the given inequality is therefore the same as
\[
\text{dim }K^{\prime}>\frac{\text{dim }V^{\prime}\text{ dim }K}{\text{dim }%
V}\text{,}%
\]
which is the unstability criterion for the pair $(K,L)$ for the $SL(V)$ action
on the space of monads. Hence the sheaf $F$ is unstable for the $SL(3)$ action
by pullback, linearized with respect to the line bundle $E_{m}$.The proof is complete.
\end{proof}

Note that when $F$ is a vector bundle there are no skyscraper sheaves in the
splitting of $F$ over a line.

Our next task is to prove that for the action of $SL(3)$ on $\mathcal{M}%
_{P^{2}}(r,0,n)$ a generic point is stable. We will prepare for this direction
by first giving examples of semistable sheaves of arbitrary high rank that are
$SL(3)$ unstable. The Serre criterion helps us construct bundles of rank $2$
with desired splitting properties when restricted to lines in $P^{2}$. Choose
$n+1$ points $x_{1},...x_{n+1}$ in $P^{2}$ such that $d+1$ of them lie on a
line $l$, $d\neq n$. Let $\mathcal{Z}$ be the zero scheme corresponding to the
$n+1$ points. Let $F^{\prime}$ be given as an extension%
\[
0\longrightarrow\mathcal{O}_{P^{2}}\longrightarrow F^{\prime}\longrightarrow
\mathcal{I}_{Z}(2)\longrightarrow0\text{.}%
\]
Let $F=F^{\prime}\otimes\mathcal{O}(-1)$. Then $c_{1}(F)=0$ and $c_{2}(F)=n$.
Also, since not all the points are on a line, $F$ is a stable bundle.
Moreover, $F|_{l}=\mathcal{O}(d)\oplus\mathcal{O}(-d)$. Choosing
$d>\frac{2n}{3}$, $d\neq n$, we have lots of examples of stable bundles that
are $SL(3)$ unstable. The following proposition enables a construction of
$SL(3)$ unstable sheaves in higher rank.

\begin{proposition}
Fix an integer $q<n$, and a line $l\subseteq P^{2}$. Then there are
$F\in\mathcal{M}_{P^{2}}(r,0,n)$ such that $F|_{l}=\oplus\mathcal{O}%
(d_{i})\oplus_{j=1,...p}C_{j}$ where $C_{j}$ are skyscraper sheaves and
\[
\Sigma_{d_{i}\geq0}d_{i}+p=q\text{.}%
\]
\end{proposition}

\begin{proof}
Note first that it is enough to prove the existence of such an $F$ for some
line $l^{\prime}\subseteq P^{2}$. For then, we can choose an automorphism $f $
of $P^{2}$ such that $f^{-1}(l^{\prime})=l$, and then $f^{\ast}(F)$ provides
the required example. Tthe proof will be by induction on the rank of the
sheaf. The statement is true for rank $2$ sheaves by the Serre construction
described above. Assume now that the statement holds for rank $r $. We need to
construct rank $r+1$ sheaves with the desired property.

Suppose that $F$ is semistable of rank $r$ such that there is a line
$l\subseteq P^{2}$ with $F|_{l}=\oplus\mathcal{O}(d_{i})\oplus_{j=1,...p}%
C_{j}$ and $\Sigma_{d_{i}\geq0}d_{i}+p=q$. Recall that $H^{1}%
(F)=Ext(\mathcal{O},F)=n-r$. Consider the extension
\[
0\longrightarrow F\longrightarrow G\longrightarrow\mathcal{O}\longrightarrow
0\text{.}
\]
Then $G$ is semistable and we have $0\rightarrow H^{0}(F(-1)|_{l})\rightarrow
H^{0}(G(-1)|_{l})\rightarrow0$. Let $G|_{l}=\oplus\mathcal{O}(e_{i}%
)\oplus_{j=1,...s}D_{j}$. Then
\[
H^{0}(G(-1)|_{l})=\Sigma_{e_{i}\geq0}e_{i}+s=H^{0}(F(-1)|_{l})=\Sigma
_{d_{i}\geq0}d_{i}+p=q\text{.}
\]
\end{proof}

Choosing $q>\frac{2n}{3}$ we can construct examples of unstable sheaves of
arbitrary rank.

\begin{theorem}
For the $SL(3)$ action on $\mathcal{M}_{P^{2}}(r,0,n)$, linearized by the line
bundle $E_{m}$ a generic point is stable.
\end{theorem}

\begin{proof}
It is enough to produce one stable point. Since the space $\mathcal{M}_{P^{2}%
}(r,0,n)$ is irreducible genericness will follow. We know, from the Hulsbergen
construction that a generic poimt in $\mathcal{M}_{P^{2}}(2,0,n)$ is stable.
The proof will be by induction on the rank of the sheaf $r\leq n$. Assume that
there is a $SL(3)$ stable point $G$ in $\mathcal{M}_{P^{2}}(r-1,0,n)$.
Consider a non split extension $0\rightarrow G\rightarrow F\rightarrow
\mathcal{O}\rightarrow0$. Then $F$ is Geiseker semistable with rank $r$ and
Chern classes $c_{1}=0$ and $c_{2}=n$ and hence determines a point in
$\mathcal{M}_{P^{2}}(r,0,n)$. We will show that $F$ is stable under the
$SL(3)$ action. In order to show this we need to prove that the pair
$(K(F),L(F)),$where $K(F)\simeq H^{1}(F(-2))$ and $\ L(F)\simeq H^{1}(F)$, is
$SL(H(F))\times SL(V)$ stable in the space of monads under the action of
$SL(H(F))\times SL(V)$, where $H(F)\simeq H^{1}(F(-1))$. Let $T_{G}%
=SL(H(G))\times SL(V)$ and $T_{F}=SL(H(F))\times SL(V)$. Note that from the
sequence defining $F$ we get $K(G)\simeq K(F)=K$(say) and $H(G)\simeq
H(F)=H$(say). Thus $T_{G}\simeq T_{F}$ $=T$ (say). For computational
convenience we will think of the space of monads as a subset of product of
grassmanians of subspaces rather than quotients. A pair $(K,L)$ corresponds to
a pair $(K,W)$ under this identification.

Let $\lambda\subseteq$ $T$ be a $1$-parameter subgroup. Then $\lambda
=\lambda_{1}.\lambda_{2}$ where $\lambda_{1}\subseteq SL(H)$, and $\lambda
_{2}\subseteq SL(V)$ are $1$-parameter subgroups. Choose $H_{i}\subseteq H$
such that $H=\oplus_{i}H_{i}$ and $\lambda_{1}(h_{i})=t^{p_{i}}h_{i}$ for
$h_{i}\in H_{i}$, and $V_{j}\subseteq V$ such that $\lambda_{2}(v_{j}%
)=t^{q_{j}}v_{j}$ for $v_{j}\in V_{j}$. Then rearragne so that
\[
H\otimes V=\oplus_{i,j}E_{i,j}=\oplus_{p=1,...,s}E_{p}%
\]
where $E_{i,j}=H_{i}\otimes V_{j}$ and $\lambda(e_{p})=t^{u_{p}}e_{p}$ for
$e_{p}\in E_{p}$ and $u_{1}>u_{2}>...>u_{s}$. Note that by construction if
$E_{p}=E_{i,j}$ then $u_{p}=p_{i}+q_{j}$. Let $F_{t}=\oplus_{p=1,...,t}E_{p}$.
Then we have a flag
\[
F_{1}\subseteq F_{2}\subseteq...\subseteq F_{s}=H\otimes V\text{.}%
\]
Define $K_{i}=K\cap F_{i}$ and $W_{i}=W\cap F_{i}$ for $\ i=1,...,s$.

We know that $\mu((K(G),W(G)),\lambda)>0$, by the stability of $G$. We need to
prove that $\mu((K(F),W(F)),\lambda)>0$. Stability of $F$ will follow from the
Hilbert-Mumford criterion. Now for the action linearized by $(p,q)$,
\[
-\mu((K,W),\lambda)=-p\mu(K(G),\lambda)-q\mu(W(G),\lambda)
\]
hence
\[
-\mu((K(G),W(G)),\lambda)=p\Sigma_{i=1,...,s}u_{i}(\text{dim}K_{i}%
-\text{dim}K_{i-1})+q\Sigma_{i=1,...,s}u_{i}(\text{dim}W_{i}(G)-\text{dim}%
W_{i-1}(G))
\]
and
\[
-\mu((K(F),W(F)),\lambda)=p\Sigma_{i=1,...,s}u_{i}(\text{dim}K_{i}%
-\text{dim}K_{i-1})+q\Sigma_{i=1,...,s}u_{i}(\text{dim}W_{i}(F)-\text{dim}%
W_{i-1}(F))\text{.}%
\]
Simplifying the equalities gives
\[
\frac{-1}{p}\mu((K(G),W(G)),\lambda)=\Sigma u_{i}(\text{dim}K_{i}%
-\text{dim}K_{i-1}+\frac{q}{p}(\text{dim}W_{i}(G)-\text{dim}W_{i-1}%
(G)))\text{.}%
\]
Recall that for the action on $\mathcal{M}_{P^{2}}(r-1,0,n)$, $p=(r-1)m-n$ and
$q=n$. For large $m$, we have $p>>q$ and hence
\[
\frac{-1}{p}\mu((K(G),W(G)),\lambda)\sim\Sigma u_{i}(dimK_{i}-dimK_{i-1}%
)<0\text{,}%
\]
since $\frac{-1}{p}\mu((K(G),W(G)),\lambda)>0$. When we consider the action on
$\mathcal{M}_{P^{2}}(r,0,n)$, $p=rm-n$ and $q=n.$Again for large $m$, $p>>q$
and we have $\frac{-1}{p}\mu((K(F),W(F)),\lambda)\sim\Sigma u_{i}($dim$K_{i}%
-$dim$K_{i-1})$. Since the right hand side has been proved to be negative, we
get $\mu((K(F),W(F)),\lambda)>0$ thus proving stability of $F.$
\end{proof}

\end{document}